\documentclass{amsart}%%NOTE: THIS IS A LATEX FILE and sent to Acta Arithmetica, 6/24/08
\usepackage{amssymb, latexsym}
\usepackage{url}

\DeclareMathOperator{\Aut}{\mathrm{Aut}}

\DeclareMathOperator{\Tor}{\mathrm{Tor}}

\begin{document}
 \bibliographystyle{plain}

 \newtheorem{theorem}{Theorem}
 \newtheorem{lemma}{Lemma}
 \newtheorem{corollary}{Corollary}
 \newtheorem{conjecture}{Conjecture}
 \newtheorem{definition}{Definition}
 \newcommand{\mc}{\mathcal}
 \newcommand{\A}{\mc{A}}
 \newcommand{\B}{\mc{B}}
 \newcommand{\cc}{\mc{C}}
 \newcommand{\D}{\mc{D}}
 \newcommand{\E}{\mc{E}}
 \newcommand{\F}{\mc{F}}
 \newcommand{\G}{\mc{G}}
 \newcommand{\FN}{\F_n}
 \newcommand{\I}{\mc{I}}
 \newcommand{\J}{\mc{J}}
 \newcommand{\eL}{\mc{L}}
 \newcommand{\M}{\mc{M}}
 \newcommand{\eN}{\mc{N}}
 \newcommand{\qq}{\mc{Q}}
 \newcommand{\U}{\mc{U}}
 \newcommand{\X}{\mc{X}}
 \newcommand{\Y}{\mc{Y}}
 \newcommand{\C}{\mathbb{C}}
 \newcommand{\R}{\mathbb{R}}
 \newcommand{\N}{\mathbb{N}}
 \newcommand{\Q}{\mathbb{Q}}
 \newcommand{\Z}{\mathbb{Z}}
 \newcommand{\ff}{\mathfrak F}
 \newcommand{\fb}{f_{\beta}}
 \newcommand{\fg}{f_{\gamma}}
 \newcommand{\gb}{g_{\beta}}
 \newcommand{\vphi}{\varphi}
 \newcommand{\bo}{\boldsymbol 0}
 \newcommand{\bgamma}{\boldsymbol \gamma}
 \newcommand{\bt}{\boldsymbol t}
 \newcommand{\bu}{\boldsymbol u}
 \newcommand{\bx}{\boldsymbol x}
 \newcommand{\bz}{\boldsymbol z}
 \newcommand{\whG}{\widehat{G}}
 \newcommand{\oK}{\overline{K}}
 \newcommand{\oKt}{\overline{K}^{\times}}
 \newcommand{\oQ}{\overline{\Q}}
 \newcommand{\oq}{\oQ^{\times}}
 \newcommand{\oQt}{\oQ^{\times}/\Tor\bigl(\oQ^{\times}\bigr)}
 \newcommand{\ot}{\Tor\bigl(\oQ^{\times}\bigr)}
 \newcommand{\dx}{\text{\rm d}x}
 \newcommand{\dy}{\text{\rm d}y}
 \newcommand{\dmu}{\text{\rm d}\mu}
 \newcommand{\dnu}{\text{\rm d}\nu}
 \newcommand{\dla}{\text{\rm d}\lambda}
 \newcommand{\dlav}{\text{\rm d}\lambda_v}
 \def\today{\number\time, \ifcase\month\or
  January\or February\or March\or April\or May\or June\or
  July\or August\or September\or October\or November\or December\fi
  \space\number\day, \number\year}

\title[Weil height]{A Banach space\\determined by the Weil height}
\author{Daniel Allcock and Jeffrey~D.~Vaaler}
\thanks{The research of both authors was supported by the National Science Foundation, 
DMS-06-00112 and DMS-06-03282, respectively.}
\address{Department of Mathematics, University of Texas, Austin, Texas
78712 USA}
\email{allcock@math.utexas.edu}
\email{vaaler@math.utexas.edu}
%\allowdisplaybreaks
\numberwithin{equation}{section}

%%%%%%%%%%%%%%%%%%%%%%%%%%%%%%%%%%%%%%%%%%%%%%%%%%%%%%%%%%%%%%%%%%%%%%%%%%%%%%%%
%%%%%%%%%%%%%%%%%%%%%%%%%%%%%%%%%%%%%%%%%%%%%%%%%%%%%%%%%%%%%%%%%%%%%%%%%%%%%%%%

\begin{abstract}  The absolute logarithmic Weil height is well defined on the
quotient group $\oQt$ and induces a
metric topology in this group.  We show that the completion of this metric space is
a Banach space over the field
$\R$ of real numbers.  We further show that this Banach space is isometrically
isomorphic to a co-dimension one subspace
of $L^1(Y, \B, \lambda)$, where $Y$ is a certain totally disconnected, locally
compact space, $\B$ is the
$\sigma$-algebra of Borel subsets of $Y$, and $\lambda$ is a certain measure
satisfying an invariance property with 
respect to the absolute Galois group $\Aut(\oQ/\Q)$.  
\end{abstract}

%%%%%%%%%%%%%%%%%%%%%%%%%%%%%%%%%%%%%%%%%%%%%%%%%%%%%%%%%%%%%%%%%%%%%%%%%%%%%%%%
%%%%%%%%%%%%%%%%%%%%%%%%%%%%%%%%%%%%%%%%%%%%%%%%%%%%%%%%%%%%%%%%%%%%%%%%%%%%%%%%

\maketitle
%\today

\section{Introduction}

Let $k$ be an algebraic number field of degree $d$ over $\Q$, $v$ a place of $k$ and
$k_v$ the completion of $k$ at $v$.  We select two absolute values from the place
$v$.  The first is denoted 
by $\|\ \|_v$ and defined as follows:
\begin{itemize}
\item[(i)] if $v|\infty$ then $\|\ \|_v$ is the unique absolute value on $k_v$ that
extends the usual absolute value on $\Q_{\infty} = \R$,
\item[(ii)] if $v|p$ then $\| \ \|_v$ is the unique absolute value on $k_v$ that
extends the usual $p$-adic absolute value on $\Q_p$.
\end{itemize}
The second absolute value is denoted by $|\ |_v$ and defined by
$|x|_v=\|x\|_v^{d_v/d}$ for all $x$ in 
$k_v$, where $d_v = [k_v:\Q_v]$ is the local degree.  If $\alpha \not = 0$ is in $k$
then these absolute values satisfy the product formula
\begin{equation}\label{prod1}
\prod_v |\alpha|_v = 1.
\end{equation}
Let $\oQ$ be an algebraic closure of $\Q$ and $\oq$ the multiplicative group of
nonzero elements in $\oQ$.
The {\it absolute, logarithmic Weil height} (or simply the {\it height}) 
\begin{equation*}\label{ht0}
h:\oq\rightarrow [0,\infty) 
\end{equation*}
is defined as follows.  Let $\alpha$ be a nonzero algebraic number, we select an algebraic number field $k$ 
containing $\alpha$, and then
\begin{equation}\label{ht1}
h(\alpha) = \sum_v \log^+|\alpha|_v,
\end{equation}
where the sum on the right of (\ref{ht1}) is over all places $v$ of $k$.
It can be shown that $h(\alpha)$ is well defined because the right hand side of (\ref{ht1})
does not depend on the field $k$.  By combining (\ref{prod1}) and (\ref{ht1}) we obtain the useful identity
\begin{equation}\label{ht2}
2 h(\alpha) = \sum_v~\bigl|\log |\alpha|_v\bigr|,
\end{equation}
where $|\ |$ (an absolute value without a subscript) is the usual archimedean absolute value on $\R$.

Let $\ot$ denote the torsion subgroup of $\oq$ and write
\begin{equation*}
\G = \oQt
\end{equation*}
for the quotient group.  If $\zeta$ is a point in $\ot$, then it is immediate from
(\ref{ht1})
that $h(\alpha) = h(\zeta\alpha)$ for all points $\alpha$ in $\oq$.  Thus $h$ is
constant on each coset of the quotient group $\G$, and so we may regard the height as a map
\begin{equation*}
h:\G\rightarrow [0,\infty).
\end{equation*}
The height has the following well known properties (see \cite[Section
1.5]{bombieri2006}):
\begin{itemize}
\item[(i)] $h(\alpha) = 0$ if and only if $\alpha$ is the identity element in $\G$,
\item[(ii)] $h(\alpha^{-1}) = h(\alpha)$ for all $\alpha$ in $\G$,
\item[(iii)] $h(\alpha\beta) \le h(\alpha) + h(\beta)$ for all $\alpha$ and $\beta$
in $\G$.
\end{itemize}
These conditions imply that the map $(\alpha, \beta)\mapsto h(\alpha \beta^{-1})$
defines a metric on the group $\G$ and therefore induces a metric topology.  Our objective in this paper is to
determine the completion of $\G$ with respect to this metric.

Let $r/s$ denote a rational number, where $r$ and $s$ are relatively prime integers
and $s$ is positive.  If $\alpha$
is in $\oq$ and $\zeta_1$ and $\zeta_2$ are in $\ot$, then all roots of the two
polynomial equations
\begin{equation*}
x^s - (\zeta_1 \alpha)^r = 0\quad\text{and}\quad x^s - (\zeta_2 \alpha)^r = 0
\end{equation*}
belong to the same coset in $\G$.  If we write $\alpha^{r/s}$ for this coset, we
find that 
\begin{equation*}
(r/s, \alpha)\mapsto \alpha^{r/s}
\end{equation*}  
defines a scalar multiplication in the abelian group $\G$.  This shows that $\G$ is a
vector space (written
multiplicatively) over the field $\Q$ of rational numbers.  Moreover, we have (see
\cite[Lemma 1.5.18]{bombieri2006})
\begin{equation}\label{norm1}
h\bigl(\alpha^{r/s}\bigr) = |r/s|h(\alpha).
\end{equation} 
Therefore the map 
$\alpha\mapsto h(\alpha)$ is a norm on the vector space $\G$ with respect to 
the usual archimedean absolute value $|\ |$ on its field $\Q$ of scalars.  From
these observations 
we conclude that the completion of $\G$ is a Banach space over the field $\R$ of
real numbers.  It 
remains now to give an explicit description of this Banach space.

Let $Y$ denote the set of all places $y$ of the field $\oQ$.  Let $k\subseteq \oQ$
be an algebraic 
number field such that $k/\Q$ is a Galois extension.  At each place $v$ of $k$ we write
\begin{equation}\label{un0}
Y(k,v) = \{y\in Y: y|v\}
\end{equation}
for the subset of places of $Y$ that lie over $v$.  Clearly we can express $Y$ as
the disjoint union
\begin{equation}\label{un1}
Y = \bigcup_v~Y(k,v),
\end{equation}
where the union is over all places $v$ of $k$.  If $y$ is a place in $Y(k,v)$ we
select an absolute value
$\|\ \|_y$ from $y$ such that the restriction of $\|\ \|_y$ to $k$ is equal to $\|\
\|_v$.  As the
restriction of $\|\ \|_v$ to $\Q$ is one of the usual absolute values on $\Q$, it
follows that this choice
of the normalized absolute value $\|\ \|_y$ does not depend on $k$.

In section 2 we show that each subset
$Y(k,v)$ can be expressed as an inverse limit of finite sets.  This determines a
totally disconnected,
compact, Hausdorff topology in $Y(k,v)$.  Then (\ref{un1}) implies that $Y$ is a
totally disconnected,
locally compact, Hausdorff space.  Again the topology in $Y$ does not depend on the
field $k$.
We also show that the absolute Galois group $\Aut(\oQ/k)$ acts transitively and
continuously on the elements 
of each compact, open subset $Y(k,v)$.

In section 4 we establish the existence of a regular measure $\lambda$, defined on
the Borel subsets of $Y$,
that is positive on open sets, finite on compact sets, and satisfies $\lambda(\tau
E) = \lambda(E)$ for all 
automorphisms $\tau$ in $\Aut(\oQ/k)$ and all Borel subsets $E$ of $Y$.  The
restriction of the measure
$\lambda$ to each subset $Y(k,v)$ is unique up to a positive multiplicative
constant.  We construct $\lambda$
so that 
\begin{equation}\label{norm2}
\lambda\bigl(Y(k,v)\bigr) = \frac{[k_v:\Q_v]}{[k:\Q]}
\end{equation}
for each Galois extension $k$ of $\Q$ and each place $v$ of $k$.  It follows from
our construction that $\lambda$ 
does not depend on the number field $k$.  In particular, if $l$ is any finite,
Galois extension of $\Q$, if 
$w$ is place of $l$ and
\begin{equation*}
Y(l,w) = \{y\in Y: y|w\},
\end{equation*}
then
\begin{equation*}
\lambda\bigl(Y(l,w)\bigr) = \frac{[l_w:\Q_w]}{[l:\Q]}.
\end{equation*}

Next we consider the real Banach space
$L^1(Y, \B, \lambda)$, where $\B$ denotes the $\sigma$-algebra of Borel subsets of
$Y$.  Let 
\begin{equation}\label{space1}
\X = \Big\{F\in L^1(Y, \B, \lambda): \int_{Y}F(y)\ \dla(y) = 0\Big\},
\end{equation}
so that $\X$ is a co-dimension one linear subspace of $L^1(Y, \B, \lambda)$.  For
each point $\alpha$ in $\G$ we
define a map $f_{\alpha}:Y\rightarrow \R$ by 
\begin{equation}\label{iso1}
f_{\alpha}(y) = \log\|\alpha\|_y.
\end{equation}
If $k$ is a finite Galois extension of $\Q$ that contains $\alpha$, then $y\mapsto
\log \|\alpha\|_y$ is
constant on each compact, open set $Y(k,v)$.  And the value of this map on each set
$Y(k,v)$ is nonzero
for only finitely many places $v$ of $k$.  It follows that $f_{\alpha}(y)$ is a
continuous function
on $Y$ with compact support.  Using (\ref{norm2}) and the product formula
(\ref{prod1}), we find that
\begin{align}\label{norm3}
\begin{split}
\int_Y f_{\alpha}(y)\ \dla(y) &= \sum_v \int_{Y(k,v)} \log \|\alpha\|_y\ \dla(y)\\
         &= \sum_v~\frac{[k_v:\Q_v]}{[k:\Q]} \log \|\alpha\|_v\\
        &= \sum_v~\log |\alpha|_v = 0.
\end{split}
\end{align}
This shows that $\alpha\mapsto f_{\alpha}(y)$ maps $\G$ into the subspace $\X$.  It
follows easily that
\begin{equation*}
f_{\alpha\beta}(y) = f_{\alpha}(y) + f_{\beta}(y)\quad\text{and}\quad
f_{\alpha^{r/s}}(y) = (r/s) f_{\alpha}(y),
\end{equation*}
and therefore the map $\alpha\mapsto f_{\alpha}(y)$ is a linear map from the vector
space $\G$ into $\X$.
The $L^1$-norm of each function $f_{\alpha}$ is given by
\begin{align}\label{norm4}
\begin{split}
\int_Y \bigl|f_{\alpha}(y)\bigr|\ \dla(y) 
        &= \sum_v \int_{Y(k,v)} \bigl|\log \|\alpha\|_y\bigr|\ \dla_v(y)\\
         &= \sum_v~\frac{[k_v:\Q_v]}{[k:\Q]}\bigl|\log \|\alpha\|_v\bigr|\\
        &= \sum_v~\bigl|\log |\alpha|_v\bigr|\\
        &= 2 h(\alpha).
\end{split}
\end{align}
This shows that the map $\alpha\mapsto f_{\alpha}$ is a linear isometry from 
the vector space $\G$ with norm determined by $2h$ into the subspace $\X$ with
the $L^1$-norm.  Let
\begin{equation}\label{norm5}
\F = \big\{f_{\alpha}(y): \alpha\in \G\big\}
\end{equation}
denote the image of $\G$ under this linear map.  Then $\alpha\rightarrow f_{\alpha}$
is a linear isometry from the vector space $\G$ (written multiplicatively) onto the 
vector space $\F$ (written additively).  Now the completion of $\G$ is determined by 
finding the closure of $\F$ in $\X$.

\begin{theorem}\label{thm1} 
Let $\X$ be the co-dimension one subspace of $L^1(Y, \B, \lambda)$ defined by {\rm
(\ref{space1})}.  Then $\F$ is dense in $\X$.
\end{theorem}

\noindent It is immediate from Theorem \ref{thm1} that  there exists
an isometric isomorphism from the completion of the vector space $\G$ with respect
to the height $2h$ onto the real Banach space $\X$.  

The functions in the vector space $\F$ belong to the real vector space $C_c(Y)$ of 
continuous functions with compact support.  Hence $\F$ belongs to the space 
$L^p(Y,\B,\lambda)$ for $1\le p \le \infty$.  Theorem \ref{thm1} asserts that the 
closure of $\F$ in $L^1(Y,\B,\lambda)$ is the
co-dimension one subspace $\X$.  We also determine the closure of $\F$ with
respect to the other $L^p$-norms. 

\begin{theorem}\label{thm2}  
If $1 < p < \infty$ then $\F$ is dense in $L^p(Y,\B,\lambda)$. 
\end{theorem}

Let $C_0(Y)$ denote the Banach space of continuous real valued functions on $Y$ which 
vanish at infinity equipped with the $\sup$-norm.  As $\F \subseteq C_c(Y) \subseteq C_0(Y)$, 
it is clear that the closure of $\F$ with respect to the $\sup$-norm is a subspace of $C_0(Y)$.

\begin{theorem}\label{thm3}
The vector space $\F$ is dense in $C_0(Y)$.
\end{theorem}

It follows from the classification of separable $L^p$-spaces 
(see \cite[pp.~14-15]{johnson2001}) that the Banach space $L^1(Y, \B, \lambda)$ has a
Schauder basis, or simply a basis.  As $\X \subseteq L^1(Y, \B, \lambda)$ is a closed subspace 
of co-dimension one, it is easy to show that $\X$ also has a basis.  Then it follows from a 
well known result of Krein, Milman and Rutman \cite{krein1940} that a basis for $\X$ can be 
selected from the dense subset $\F$.  Thus there exists a sequence of distinct elements
$\alpha_1, \alpha_2, \dots $ in $\G$ such that the corresponding collection of functions
\begin{equation}\label{basis1}
\big\{f_{\alpha_1}(y), f_{\alpha_2}(y), \dots \big\}
\end{equation}
is a basis for the Banach space $\X$.  That is, for every function $F$ in $\X$ there exists a
{\it unique} sequence of real numbers $x_1, x_2, \dots $ such that
\begin{equation*}\label{basis2}
F(y) = \lim_{N\rightarrow \infty} \sum_{n=1}^N x_n f_{\alpha_n}(y)
\end{equation*}  
in $L^1$-norm.  While these remarks establish the existence of such a basis, it would be of interest 
to construct an explicit example of a sequence $\alpha_1, \alpha_2, \dots $ in $\G$ such that the 
the corresponding sequence of functions (\ref{basis1}) forms a basis for $\X$.

%%%%%%%%%%%%%%%%%%%%%%%%%%%%%%%%%%%%%%%%%%%%%%%%%%%%%%%%%%%%%%%%%%%%%%%%%%%%%%%%%%%%%%%
%%%%%%%%%%%%%%%%%%%%%%%%%%%%%%%%%%%%%%%%%%%%%%%%%%%%%%%%%%%%%%%%%%%%%%%%%%%%%%%%%%%%%%%

\section{Preliminary Lemmas}

We have stated Theorem \ref{thm1} for the Weil height on algebraic number fields. 
However, many of the arguments
can be given in the more general setting of a field $K$ with a proper set of
absolute values satisfying a 
product formula.  We now describe this situation.

Let $K$ be a field and let $v$ be a place of $K$.  That is, $v$ is an equivalence
class of nontrivial absolute
values on $K$.  We write $K_v$ for the completion of $K$ at the place $v$.  If $L/K$
is a finite extension of
fields then there exist finitely many places $w$ of $L$ such that $w|v$.  In general
we have
\begin{equation*}
\sum_{w|v}~[L_w:K_v] \le [L:K],
\end{equation*}
where $L_w$ is the completion of $L$ at $w$.  We say that $v$ is {\it well behaved}
if the identity
\begin{equation*}
\sum_{w|v}~[L_w:K_v] = [L:K]
\end{equation*}
holds for all finite extensions $L/K$, (see \cite[Chapter 1, section 4]{lang1983}). 

Let $\M_K$ be a collection of distinct places of $K$ and at each place $v$ in $\M_K$
let $\|\ \|_v$ denote an absolute
value from $v$.  We say that the collection of absolute values
\begin{equation}\label{prod2}
\big\{\|\ \|_v: v\in\M_K\big\}
\end{equation}
is {\it proper} if it satisfies the following conditions:
\begin{itemize}
\item[(i)] each place $v$ in $\M_K$ is well behaved, 
\item[(ii)] if $\alpha$ is in $K^{\times}$ then $\|\alpha\|_v \not= 1$ for at most
finitely many places $v$ in $\M_K$,
\item[(iii)]  if $\alpha$ is in $K^{\times}$ then the absolute values in
(\ref{prod2}) satisfy the product formula
\begin{equation*}
\prod_{v\in\M_K} \|\alpha\|_v = 1.
\end{equation*}
\end{itemize}

Now suppose that (\ref{prod2}) is a proper set of absolute values on $K$ and $L/K$
is a finite extension of fields.
Let $\M_L$ be the collection of places of $L$ that extend the places in $\M_K$. 
That is, if $W_v(L/K)$ is the
finite set of places $w$ of $L$ such that $w|v$, then
\begin{equation*}
\M_L = \bigcup_{v\in\M_K} W_v(L/K).
\end{equation*}
At each place $w$ in $W_v(L/K)$ we select an absolute value $\|\ \|_w$ that extends
the absolute value $\|\ \|_v$ 
on $K$.  Then we define an equivalent absolute value $|\ |_w$ from the place $w$ by
setting
\begin{equation*}
\log |\alpha|_w = \frac{[L_w:K_v]}{[L:K]} \log \|\alpha\|_w
\end{equation*}
for all $\alpha$ in $L^{\times}$.  In general $\|\ \|_w$
 and $|\ |_w$ are distinct but equivalent absolute values on $L$.
And we note that $|\ |_w$ {\it is} an absolute value because
\begin{equation*}
0 < \frac{[L_w:K_v]}{[L:K]} \le 1.
\end{equation*}
Then it follows, as in \cite[Chapter 2, section 1]{lang1983}, that
\begin{equation}\label{prod3}
\big\{|\ |_w: w\in\M_L\big\}
\end{equation}
is a proper set of absolute values on $L$.  In particular, if $\alpha$ is in
$L^{\times}$ then the absolute values
in (\ref{prod3}) satisfy the product formula
\begin{equation*}\label{prod4}
\prod_{w\in M_L} |\alpha|_w = 1.
\end{equation*}

%If $\oK$ is an algebraic closure of $K$ and $\alpha$ is in $\oKt$, we define
%the absolute, logarithmic height of $\alpha$ with respect to $\M_K$ by
%\begin{equation}\label{prod5}
%h(\alpha) = \sum_{w\in\M_L} \log^+ |\alpha|_w,
%\end{equation}
%where $L/K$ is a finite extension that contains $\alpha$.  It can be shown, as in
%\cite[Chapter 3, section 1]{lang1983}, 
%that the sum on the right of (\ref{prod5}) does not depend on $L$, and therefore
%$h:\oKt\rightarrow [0,\infty)$ is 
%well defined.

We assume that $K\subseteq N$ are fields, that $N/K$ is a (possibly infinite) Galois
extension, 
and we write $\Aut(N/K)$ for the corresponding Galois group.  We give
$\Aut(N/K)$ the Krull topology, and we briefly recall how this is defined.  Let
$\eL$ denote the 
set of intermediate fields $L$ such that $K\subseteq L\subseteq N$ and $L/K$ is a
finite Galois extension.  Obviously 
$\eL$ is partially ordered by set inclusion.  If $L$ and $M$ are in $\eL$ then the
composite field $LM$ is in 
$\eL$, $L\subseteq LM$, $M\subseteq LM$, and therefore $\eL$ is a directed set.  For
each $L$ in $\eL$ let 
$\Aut(L/K)$ denote the Galois group of automorphisms of $L$ that fix $K$.  If
$L\subseteq M$ are both in 
$\eL$, we define $\pi_L^M:\Aut(M/K)\rightarrow \Aut(L/K)$ to be the map that
restricts the domain of an automorphism 
in $\Aut(M/K)$ to the subfield $L$.  Then each map $\pi_L^M$ is a surjective
homomorphism of groups and $\pi_L^L$ is the 
identity map.  It follows that
\begin{equation*}
\big\{\Aut(L/K), \pi_L^M\big\}
\end{equation*}
is an inverse system, and $\Aut(N/K)$ can be identified with the inverse (or
projective) limit:
\begin{equation*}
\Aut (N/K) = \lim_{\substack{\longleftarrow\\L\in\eL}} \Aut(L/K).
\end{equation*}
Thus $\Aut(N/K)$ is a profinite group, and therefore is a totally disconnected,
compact, Hausdorff, topological 
group.  We write 
\begin{equation*}
\pi_L:\Aut(N/K)\rightarrow\Aut(L/K) 
\end{equation*}
for the canonical map associated with each $L$ in $\eL$.  Then $\pi_L$ is continuous
and the collection of open sets
\begin{equation}\label{base1}
\big\{\pi_L^{-1}(\tau): L\in\eL\ \text{and}\ \tau\in\Aut(L/K)\big\}
\end{equation}
is a basis for the Krull topology in $\Aut(N/K)$.

Next we assume that $v$ is a place of the field $K$.  That is, $v$ is an equivalence
class of nontrivial absolute 
values on $K$.  If $L$ is in $\eL$ we write $W_v(L/K)$ for the set of places $w$ of
$L$ such 
that $w|v$.  As $L/K$ is a finite extension, it follows that $W_v(L/K)$ is a finite
set.  If $L\subseteq M$
belong to $\eL$ we define connecting maps 
\begin{equation*}
\psi_L^M:W_v(M/K)\rightarrow W_v(L/K) 
\end{equation*}
as follows: if $w_M$ belongs 
to $W_v(M/K)$ then $\psi_L^M(w_M)$ is the unique place $w_L$ in $W_v(L/K)$ such that
$w_M|w_L$.  If $L\subseteq M$ 
are in $\eL$ then each absolute value on $L$ extends to $M$ and therefore each
connecting map $\psi_L^M$ 
is surjective.  We give each finite set $W_v(L/K)$ the discrete topology so that
each map $\psi_L^M$ is continuous. 
Clearly $\psi_L^L$ is the identity map.  We find that
\begin{equation*}
\big\{W_v(L/K), \psi_L^M\big\}
\end{equation*}
is an inverse system of finite sets.  Let
\begin{equation*}
Y(K,v) = \lim_{\substack{\longleftarrow\\L\in\eL}} W_v(L/K)
\end{equation*}   
denote the inverse limit and write $\psi_L:Y(K,v)\rightarrow W_v(L/K)$
for the canonical continuous map associated to each $L$ in $\eL$.  It follows, as 
in \cite[Appendix 2, section 2.4]{dugundji1968}, that $Y(K,v)$ is a nonempty,
totally disconnected, compact, 
Hausdorff space.  Moreover, see \cite[Appendix 2, section 2.3]{dugundji1968}, the
collection of open sets
\begin{equation}\label{base2}
\big\{\psi_L^{-1}(w): L\in\eL\ \text{and}\ w\in W_v(L/K)\big\}
\end{equation}
is a basis for the topology of $Y(K,v)$.  Clearly each subset in the collection
(\ref{base2}) is also compact,
and for each field $L$ in $\eL$ we can write
\begin{equation*}\label{base3}
Y(K,v) = \bigcup_{w\in W_v(L/K)} \psi_L^{-1}(w)
\end{equation*}
as a disjoint union of open and compact sets.

We recall that a map $g:Y(K,v)\rightarrow \R$ is {\it locally constant} if at each
point $y$ in $Y(K,v)$ there
exists an open neighborhood of $y$ on which $g$ is constant.  

\begin{lemma}\label{lem1}  Let $g:Y(K,v)\rightarrow \R$ be locally constant.  Then
there exists $L$ in 
$\eL$ such that for each place $w$ in $W_v(L/K)$ the function $g$ is constant on the
set $\psi_L^{-1}(w)$.
\end{lemma}

\begin{proof}  At each point $y$ in $Y(K,v)$ there exists a field $L^{(y)}$ in $\eL$
and a place $w^{(y)}$
in $W_v\bigl(L^{(y)}/K\bigr)$ such that $y$ is contained in
$\psi_{L^{(y)}}^{-1}\bigl(w^{(y)}\bigr)$ and 
$g$ is constant on the open set
$\psi_{L^{(y)}}^{-1}\bigl(w^{(y)}\bigr)$.  By compactness there exists a finite
collection of fields 
$L^{(1)}, L^{(2)}, \dots , L^{(J)}$ in $\eL$, and for each integer $j$ a
corresponding place $w^{(j)}$ in 
$W_v\bigl(L^{(j)}/K\bigr)$, such that 
\begin{equation*}\label{base4}
Y(K,v) \subseteq \bigcup_{j=1}^J \psi_{L^{(j)}}^{-1}\bigl(w^{(j)}\bigr),
\end{equation*}
and $g$ is constant on each open set $\psi_{L^{(j)}}^{-1}\bigl(w^{(j)}\bigr)$.
Let $L = L^{(1)}L^{(2)} \cdots L^{(J)}$ be the composite field, which is obviously
in $\eL$.  If $w$ is a place
of $L$ then there exists an integer $j$ such that 
\begin{equation*}\label{base5}
\psi_L^{-1}(w)\cap \psi_{L^{(j)}}^{-1}\bigl(w^{(j)}\bigr)
\end{equation*}
is not empty.  As $L$ is a finite extension of $L^{(j)}$, we conclude that
$\psi_{L^{(j)}}^L(w) = w^{(j)}$, and 
therefore
\begin{equation}\label{base6}
\psi_L^{-1}(w)\subseteq \psi_{L^{(j)}}^{-1}\bigl(w^{(j)}\bigr).
\end{equation}
Then (\ref{base6}) implies that $g$ is constant on $\psi_L^{-1}(w)$.
\end{proof}

Let $C\bigl(Y(K,v)\bigr)$ denote the real Banach algebra of continuous functions 
\begin{equation*}\label{base7}
F:Y(K,v)\rightarrow \R
\end{equation*}
with the supremum norm.   Let $LC\bigl(Y(K,v)\bigr)\subseteq C\bigl(Y(K,v)\bigr)$
denote the subset of
locally constant functions.  

\begin{lemma}\label{lem2}  The subset $LC\bigl(Y(K,v)\bigr)$ is a dense subalgebra
of $C\bigl(Y(K,v)\bigr)$.
\end{lemma}

\begin{proof}  It is obvious that $LC\bigl(Y(K,v)\bigr)$ is a subalgebra of
$C\bigl(Y(K,v)\bigr)$, and that
$LC\bigl(Y(K,v)\bigr)$ contains the constant functions.  Now suppose that
$y_1$ and $y_2$ are distinct points in $Y(K,v)$.  Let $U_1$ be an open neighborhood
of $y_1$, and $U_2$ an open
neighborhood of $y_2$, such that $U_1$ and $U_2$ are disjoint.  Then there exists a
field $L$ in $\eL$ and 
a place $w$ in $W_v(L/K)$ such that 
\begin{equation*}\label{base8}
y_1\in \psi_L^{-1}(w),\quad\text{and}\quad \psi_L^{-1}(w)\subseteq U_1.
\end{equation*}
As $\psi_L^{-1}(w)$ is both open and compact,
the characteristic function of the set $\psi_L^{-1}(w)$ is a locally constant function
that separates the points $y_1$ and $y_2$.  Then it follows from the
Stone-Weierstrass theorem that the 
subalgebra $LC\bigl(Y(K,v)\bigr)$ is dense in $C\bigl(Y(K,v)\bigr)$.
\end{proof}

We select an absolute value from the place $v$ of $K$ and denote it by $\|\ \|_v$. 
If $L$ is in $\eL$ and $w$ is a
place in $W_v(L/K)$, we select an absolute value $\|\ \|_w$ from $w$ such that the
restriction of $\|\ \|_w$ to
$K$ is equal to $\|\ \|_v$.  As
\begin{equation*}
N = \bigcup_{L\in\eL}~L,
\end{equation*}
it follows that each point $(w_L)$ in $Y(K,v)$ determines a unique absolute value on
the field $N$.  That is, each
point $(w_L)$ in $Y(K,v)$ determines a unique place $y$ of $N$ such that $y|v$.  

Now suppose $y$ is a place of $N$
such that $y|v$.  Select an absolute value $\|\ \|_y$ from $y$ such that the
restriction of $\|\ \|_y$ to the 
subfield $K$ is equal to $\|\ \|_v$.  If $L$ is in $\eL$ then the restriction of
$\|\ \|_y$ to $L$ must equal 
$\|\ \|_{w_L}$ for a unique place $w_L$ in $W_v(L/K)$.  Thus each place $y$ of $N$
with $y|v$ determines a
unique point $(w_L)$ in the product
\begin{equation*}
\prod_{L\in\eL} W_v(L/K)
\end{equation*} 
such that $y|w_L$ for each $L$.  It is trivial to check that 
\begin{equation*}
\psi_L^M(w_M) = w_L
\end{equation*}
whenever $L\subseteq M$ are in $\eL$.  Therefore each place $y$ of $N$ with $y|v$
determines a unique point
$(w_L)$ in the inverse limit $Y(K,v)$.  In view of these remarks we may identify
$Y(K,v)$ with the set of all 
places $y$ of $N$ that lie over the place $v$ of $K$.  In this way we determine a
totally disconnected, compact, 
Hausdorff, topology in the set of all places $y$ of $N$ that lie over the place $v$
of $K$.

%%%%%%%%%%%%%%%%%%%%%%%%%%%%%%%%%%%%%%%%%%%%%%%%%%%%%%%%%%%%%%%%%%%%%%%%%%%%%%%
%%%%%%%%%%%%%%%%%%%%%%%%%%%%%%%%%%%%%%%%%%%%%%%%%%%%%%%%%%%%%%%%%%%%%%%%%%%%%%%

\section{Galois action on places}

Next we recall that the Galois group $\Aut(N/K)$ acts on the set $Y(K,v)$ of all
places of $N$ that lie over 
$v$.  More precisely, if $\tau$ is in $\Aut(N/K)$ and $y$ is in $Y(K,v)$, then the map
\begin{equation}\label{act1}
\alpha\mapsto \|\tau^{-1}\alpha\|_y
\end{equation}
is an absolute value on $N$, and the restriction of this absolute value to $K$ is
clearly equal to $\|\ \|_v$.  
Therefore (\ref{act1}) determines a unique place $\tau y$ in $Y(K,v)$.  That is, the
identity
\begin{equation}\label{act2}
\|\tau^{-1}\alpha\|_y = \|\alpha\|_{\tau y}
\end{equation}
holds for all $\alpha$
 in $N$, for all $\tau$ in $\Aut(N/K)$, and for all places $y$ in $Y(K,v)$.  It is 
immediate that $1y = y$ and $(\sigma\tau)y = \sigma(\tau y)$ for all $\sigma$ and
$\tau$ in $\Aut(N/K)$.  Thus 
$(\tau, y)\mapsto \tau y$ defines an action of the group $\Aut(N/K)$ on the set
$Y(K,v)$.  Moreover, $\Aut(N/K)$ 
acts transitively on $Y(K,v)$, (see \cite[Chapter II, Proposition
9.1]{neukirch1999}.)  

\begin{lemma}\label{lem3}
The function $(\tau, y)\mapsto \tau y$ from $\Aut(N/K)\times Y(K,v)$ onto $Y(K,v)$
is continuous.
\end{lemma}

\begin{proof}  Let $L$ be in $\eL$ and $w$ in $W_v(L/K)$.  In view of (\ref{base2})
we must show that
\begin{equation*}
\big\{(\tau, y)\in \Aut(N/K)\times Y(K,v): \tau y\in \psi_L^{-1}(w)\big\}
\end{equation*}
is open in $\Aut(N/K)\times Y(K,v)$ with the product topology.  For $w$ in
$W_v(L/K)$ we define
\begin{equation*}
E_w = \big\{(\sigma, z)\in \Aut(L/K)\times W_v(L/K): \sigma z = w\big\}.
\end{equation*}
Then we have
\begin{align*}
\big\{(\tau, y)&\in \Aut(N/K)\times Y(K,v): \tau y\in \psi_L^{-1}(w)\big\}\\
        &= \big\{(\tau, y)\in \Aut(N/K)\times Y(K,v): \pi_L(\tau)\psi_L(y) = w\big\}\\
        &= \bigcup_{(\sigma, z)\in E_w}\big\{(\tau, y)\in \Aut(K/k)\times Y(K,v): \pi_L(\tau) 
                = \sigma\ \text{and}\ \psi_L(y) = z\big\}\\
        &= \bigcup_{(\sigma, z)\in E_w} \pi_L^{-1}(\sigma)\times \psi_L^{-1}(z),
\end{align*}
which is obviously an open subset of $\Aut(N/K)\times Y(K,v)$.
\end{proof}

%%%%%%%%%%%%%%%%%%%%%%%%%%%%%%%%%%%%%%%%%%%%%%%%%%%%%%%%%%%%%%%%%%%%%%%%%%%%%%%%
%%%%%%%%%%%%%%%%%%%%%%%%%%%%%%%%%%%%%%%%%%%%%%%%%%%%%%%%%%%%%%%%%%%%%%%%%%%%%%%%

\section{The invariant measure}

In this section it will be convenient to write $G = \Aut(N/K)$.
Let $\mu$ denote a Haar measure on the Borel subsets of the compact topological
group $G$ normalized so that 
$\mu(G) = 1$.  If $F$ is in $C\bigl(Y(K,v)\bigr)$ and $z_1$ is a point in $Y(K,v)$
then it follows 
from Lemma \ref{lem3} that $\tau\mapsto F(\tau z_1)$ is a continuous function on $G$
with values in $\R$.  Let 
$z_2$ be a second point in $Y(K,v)$.  Because $G$ acts transitively on $Y(K,v)$,
there exists $\eta$ in $G$ 
so that $\eta z_2 = z_1$.  Then using the translation invariance of Haar measure we get
\begin{equation}\label{int1}
\int_G F(\tau z_1)\ \dmu(\tau) = \int_G F(\tau\eta z_2)\ \dmu(\tau) = \int_G F(\tau
z_2)\ \dmu(\tau).
\end{equation}
It follows that the map $I_v:C\bigl(Y(K,v)\bigr)\rightarrow \R$ given by
\begin{equation}\label{int2}
I_v(F) = \int_G F(\tau z_v)\ \dmu(\tau)
\end{equation}
does not depend on the point $z_v$ in $Y(K,v)$. 

Let $\M_K$ be a collection of distinct places of $K$ and at each place $v$ in $\M_K$
let $\|\ \|_v$ denote an
absolute value from $v$.  We assume that
\begin{equation*}\label{ms1}
\big\{\|\ \|_v: v\in\M_K\big\}
\end{equation*}
is a proper collection of absolute values.  Again we assume that $N/K$ is a
(possibly infinite) Galois extension
of fields.  Let $Y$ be defined by the disjoint union
\begin{equation}\label{ms2}
Y = \bigcup_{v\in \M_K} Y(K,v).
\end{equation}
Thus $Y$ is the collection of all places $y$ of $N$ such that $y|v$ for some place
$v$ in $\M_K$.
It follows that $Y$ is a nonempty, totally disconnected, locally compact, Hausdorff
space. 

Let $C_c(Y)$ denote the real vector space of continuous functions $F:Y\rightarrow
\R$ having compact support.  If $F$ 
belongs to $C_c(Y)$ then there exists a finite subset $S_F\subseteq \M_K$ such that
$F$ is supported on the compact set
\begin{equation*}\label{ms3}
\bigcup_{v\in S_F} Y(K,v).
\end{equation*}
In particular we have $I_v(F) = 0$ for almost all places $v$ of $\M_K$.  Therefore
we define $I:C_c(Y)\rightarrow \R$ by
\begin{equation}\label{ms4}
I(F) = \sum_{v\in \M_K} \int_G F(\tau z_v)\ \dmu(\tau),
\end{equation}
where $z_v$ is a point in $Y(K,v)$ for each place $v$ in $\M_K$.  By our previous
remarks the value of each
integral on the right of (\ref{ms4}) does not depend on $z_v$, and only finitely
many integrals on the right of (\ref{ms4}) are nonzero.  Hence there is no question
of convergence in the
sum on the right of (\ref{ms4}).  

\begin{theorem}\label{thm4}
There exists a $\sigma$-algebra $\Y$ of subsets of $Y$, that contains the $\sigma$-algebra 
$\B$ of Borel sets in $Y$, and a unique, regular measure $\lambda$ defined on $\Y$, such that 
\begin{equation}\label{rr1}
I(F) = \int_Y F(y)\ \dla(y)
\end{equation}
for all $F$ in $C_c(Y)$.  Moreover, the measure $\lambda$ satisfies the following
conditions:
\begin{itemize}
\item[(i)] If $\eta$ is in $G$ and $F$ is in $L^1(Y,\Y,\lambda)$ then
\begin{equation}\label{rr2}
\int_{Y(K,v)} F(\eta y)\ \dla(y) = \int_{Y(K,v)} F(y)\ \dla(y)
\end{equation}
at each place $v$ in $\M_K$.
\item[(ii)] If $E$ is in $\Y$ then
\begin{equation*}
\lambda(E) = \inf\{\lambda(U): E\subseteq U\subseteq Y\ \text{and}\ U\ \text{is
open}\}.
\end{equation*} 
\item[(iii)] If $E$ is in $\Y$ then
\begin{equation*}
\lambda(E) = \sup\{\lambda(V): V\subseteq E\ \text{and}\ V\ \text{is compact}\}.
\end{equation*}
\item[(iv)] If $E$ is in $\Y$ and $\lambda(E) = 0$ then every subset of $E$ is in $\Y$.
\end{itemize}
\end{theorem}

\begin{proof}  Clearly (\ref{ms4}) defines a positive linear functional on $C_c(Y)$.  
By the Riesz representation theorem, (see \cite[Theorem 2.14 and Theorem 2.17]{rudin1987}),
there exists a $\sigma$-algebra $\Y$ of subsets of $Y$, containing 
the $\sigma$-algebra $\B$ of Borel sets in $Y$, and a regular measure $\lambda$ defined on
$\Y$, such that
\begin{equation}\label{int3}
I(F) = \int_Y F(y)\ \dla(y)
\end{equation}
for all $F$ in $C_c(Y)$.  If $\eta$ is in $G$ and $F$ is in $C_c(Y)$, then by the
translation invariance 
of the Haar measure $\mu$ we have
\begin{align}\label{rr3}
\begin{split}
\int_{Y(K,v)} F(\eta y)\ \dla(y) &= \int_G F(\eta\tau z)\ \dmu(\tau)\\
                          &= \int_G F(\tau z)\ \dmu(\tau)\\
                          &= \int_{Y(K,v)} F(y)\ \dla(y)
\end{split}
\end{align}
at each place $v$ in $\M_K$.
Initially (\ref{rr3}) holds for all functions $F$ in $C_c(Y)$.  As $C_c(Y)$ is 
dense in the $L^1(Y,\Y,\lambda)$, (see \cite[Theorem 3.14]{rudin1987}), it follows
in a standard manner 
that (\ref{rr3}) holds also for functions $F$ in $L^1(Y,\Y,\lambda)$. 

The properties (ii), (iii) and (iv) attributed to $\lambda$ are all consequences of 
the Riesz theorem.
\end{proof}  

Because the Haar measure $\mu$ satisfies $\mu(G) = 1$, it is immediate from
(\ref{int2}) and (\ref{rr1}) that
$\lambda\bigl(Y(K,v)\bigr) = 1$ at each place $v$ in $\M_K$.  As the places in
$\M_K$ are well behaved, 
we obtain a further identity for the $\lambda$-measure of basic open sets in each
subset $Y(K,v)$. 

\begin{theorem}\label{thm5}  If $L$ is in $\eL$ and $w$ is a place in $W_v(L/K)$, then
\begin{equation}\label{rr4}
\lambda\bigl(\psi_L^{-1}(w)\bigr) = \frac{[L_w:K_v]}{[L:K]}.
\end{equation} 
\end{theorem}

\begin{proof}  Let $\tau$ be in $G$.  Then we have
\begin{align}\label{rr5}
\begin{split}
\tau \psi_L^{-1}(w) &= \big\{\tau y \in Y(K,v): \psi_L(y) = w\big\}\\
        &= \big\{y\in Y(K,v): \pi_L\bigl(\tau^{-1}\bigr)\psi_L(y) = w\big\}\\
        &= \big\{y\in Y(K,v): \psi_L(y) = \pi_L(\tau) w\big\}\\
        &= \psi_L^{-1}\bigl(\pi_L(\tau) w\bigr).
\end{split}
\end{align} 
Now let $w_1$ and $w_2$ be distinct places in $W_v(L/K)$.  Select $\tau$ in $G$ so
that $\pi_L(\tau) w_2 = w_1$.
Then (\ref{rr5}) implies that
\begin{equation*}
\tau \psi_L^{-1}(w_2) = \psi_L^{-1}(w_1),
\end{equation*}
and using (\ref{rr2}) we find that
\begin{equation*}
\lambda\big\{\psi_L^{-1}(w_2)\big\} = \lambda\big\{\psi_L^{-1}(w_1)\big\}.
\end{equation*}
Because
\begin{equation}\label{rr6}
Y(K,v) = \bigcup_{w\in W_v(L/K)} \psi_L^{-1}(w) 
\end{equation}
is a disjoint union of $|W_v(L/K)|$ distinct sets, the sets on the right of
(\ref{rr6}) all have equal 
$\lambda$-measure, and $\lambda(Y(K,v)) = 1$, we conclude that
\begin{equation}\label{rr7}
\lambda\bigl(\psi_L^{-1}(w)\bigr) = |W_v(L/K)|^{-1}.
\end{equation}
As $v$ is well behaved we have
\begin{equation}\label{rr8}
[L:K] = \sum_{w\in W_v(L/K)} [L_w:K_v].
\end{equation}
Because $L/K$ is a Galois extension, all local degrees $[L_w:K_v]$ for $w$ in
$W_v(L/K)$ are equal, and we conclude 
from (\ref{rr8}) that
\begin{equation}\label{rr9}
|W_v(L/K)| = \frac{[L:K]}{[L_w:K_v]}.
\end{equation}
The identity (\ref{rr4}) follows now from (\ref{rr7}) and (\ref{rr9}).
\end{proof}

Let $LC_c(Y)$ be the algebra of locally constant, real valued functions on $Y$
having compact support.  Clearly
we have $LC_c(Y)\subseteq C_c(Y)$.

\begin{lemma}\label{lem4}  Let $g$ belong to $LC_c(Y)$.  Then there exists $L$ in
$\eL$ such that for each place 
$w$ in $\M_L$ the function $g$ is constant on the set $\psi_L^{-1}(w)$.
\end{lemma}

\begin{proof}  Let $S_g\subset \M_K$ be a finite set of places of $K$ such that the
support of $g$ is contained
in the compact set
\begin{equation*}
V_g = \bigcup_{v\in S_g} Y(K,v).
\end{equation*}
For each place $v$ in $S_g$ we apply Lemma \ref{lem1} to the restriction of $g$ to
$Y(K,v)$.  Thus there exists 
a field $L^{(v)}$ in $\eL$ such that for each place $w^{\prime}$ in
$W_v\bigl(L^{(v)}/K\bigr)$, the function $g$ is
constant on $\psi_{L^{(v)}}^{-1}\bigl(w^{\prime}\bigr)$.  Let $L$ be the compositum
of the finite collection of fields 
\begin{equation*}\label{rr10}
\big\{L^{(v)}: v\in S_g\big\}.
\end{equation*}
Clearly $L$ belongs to $\eL$.  

Let $w$ be a place in $\M_L$.  If $w|v$ and $v\notin S_g$, then $g$ is identically
zero on $\psi_L^{-1}(w)$,
and in particular it is constant on this set.  If $w|v$ and $v\in S_g$, then
$w|w^{\prime}$ for a unique place
$w^{\prime}$ in $W_v\bigl(L^{(v)}/K\bigr)$.  Because
\begin{equation*}\label{rr11}
\psi_L^{-1}(w) \subseteq \psi_{L^{(v)}}^{-1}\bigl(w^{\prime}\bigr)
\end{equation*} 
and $g$ is constant on $\psi_{L^{(v)}}^{-1}\bigl(w^{\prime}\bigr)$, it is obvious
that $g$ is constant on 
$\psi_L^{-1}(w)$.
\end{proof}

\begin{lemma}\label{lem5}  For $1\le p < \infty$ the set $LC_c(Y)$ is dense in
$L^p(Y,\B,\lambda)$.
And $LC_c(Y)$ is dense in $C_0(Y)$ with respect to the $\sup$-norm.
\end{lemma}

\begin{proof}
Let $1\le p < \infty$.  Because $C_c(Y)$ is dense in $L^p(Y,\B,\lambda)$, it suffices
to show that if 
$F$ is in $C_c(Y)$ and $\epsilon > 0$, then there exists a function $g$ in $LC_c(Y)$
such that
\begin{equation*}\label{rr20}
\bigg\{\int_Y |F(y) - g(y)|^p\ \dla(y)\bigg\}^{1/p} < \epsilon.
\end{equation*}
Let $S_F\subseteq \M_K$ be a nonempty, finite set of places such that $F$ is
supported on the compact set
\begin{equation*}\label{rr21}
V_F = \bigcup_{v\in S_F} Y(K,v).
\end{equation*}
For each $v$ in $S_F$ we apply Lemma \ref{lem2} to the restriction of $F$ to
$Y(K,v)$.  Thus there exists
a locally constant function $g_v:Y(K,v)\rightarrow \R$ such that
\begin{equation}\label{rr22}
\sup\big\{|F(y) - g_v(y)|: y\in Y(K,v)\big\} < |S_F|^{-1/p}\epsilon.
\end{equation} 
Now define $g:Y\rightarrow \R$ by
\begin{equation}\label{rr23}
g(y) = \begin{cases} g_v(y)& \text{if $y\in Y(K,v)$ and $v\in S_F$,}\\
                        0& \text{if $y\in Y(K,v)$ and $v\notin S_F$.}\end{cases}
\end{equation}
Then $g$ is locally constant and supported on the compact set $V_F$.  Therefore $g$
belongs to $LC_c(Y)$.
As $\lambda\bigl(Y(K,v)\bigr) = 1$ at each place $v$ in $\M_K$, we get
\begin{align*}\label{rr24}
\begin{split}
\bigg\{\int_Y |F(y) - g(y)|^p\ \dla(y)\bigg\}^{1/p} 
        &= \bigg\{\sum_{v\in S_F} \int_{Y(K,v)} |F(y) - g_v(y)|^p\ \dla(y)\bigg\}^{1/p} \\
        &< \bigg\{\sum_{v\in S_F} |S_F|^{-1} \epsilon^p\bigg\}^{1/p} \le \epsilon.
\end{split}
\end{align*}
This proves the first assertion of the lemma.  

As $C_c(Y)$ is dense in $C_0(Y)$ with respect to the $\sup$-norm, the second
assertion of the lemma follows by the
same argument.  In this case we select the locally constant functions
$g_v:Y(K,v)\rightarrow \R$ so that
\begin{equation*}\label{rr25}
\sup\big\{|F(y) - g_v(y)|: y\in Y(K,v)\big\} < \epsilon.
\end{equation*}
Then we define $g:Y\rightarrow \R$ as in (\ref{rr23}).  Again we find that $g$
belongs to $LC_c(Y)$, and the inequality
\begin{equation*}
\sup\big\{|F(y) - g(y)|: y\in Y\big\} < \epsilon
\end{equation*}
is obvious.
\end{proof}

%%%%%%%%%%%%%%%%%%%%%%%%%%%%%%%%%%%%%%%%%%%%%%%%%%%%%%%%%%%%%%%%%%%%%%%%%%%%%%%%%%
%%%%%%%%%%%%%%%%%%%%%%%%%%%%%%%%%%%%%%%%%%%%%%%%%%%%%%%%%%%%%%%%%%%%%%%%%%%%%%%%%%

\section{The completion of $\G$}

In this section we return to the situation considered in the introduction.
We let $K = \Q$, $N = \oQ$, and we let $\M_{\Q}$ be the set of all places of $\Q$.  
Then $Y$ is the set of all places of $\oQ$, and $Y$ is a nonempty, totally
disconnected, locally compact, Hausdorff space.  By Theorem \ref{thm4} there exists a
$\sigma$-algebra $\Y$ of subsets of $Y$, containing the $\sigma$-algebra $\B$ of Borel sets in $Y$,
and a measure $\lambda$ on $\Y$, satisfying the conclusions of that result.  The
basic identity (\ref{norm2}) is verified by Theorem \ref{thm5}.  Then the map
\begin{equation}\label{com1}
\alpha\rightarrow f_{\alpha}(y)
\end{equation} 
defined by (\ref{iso1}) is a linear map from the $\Q$-vector space 
\begin{equation*}\label{com2}
\G = \oQt
\end{equation*}
(written multiplicatively) into the vector space $C_c(Y)$.  The identity
(\ref{norm3}) implies that each function 
$f_{\alpha}(y)$ belongs to the closed subspace $\X\subseteq L^1(Y, \B, \lambda)$
defined by (\ref{space1}).  
It follows from basic properties of the height, and in particular (\ref{norm1}), that 
\begin{equation*}\label{com3}
\alpha\rightarrow 2h(\alpha),
\end{equation*} 
defines a norm on $\G$ with respect to the usual archimedean absolute value on $\Q$.
Then (\ref{norm4}) shows that (\ref{com1}) defines a linear isometry of $\G$ into the 
subspace $\X$.  

\begin{lemma}\label{lem6}  Let $k$ be an algebraic number field and let $v\rightarrow t_v$
be a real valued function defined on the set of all places $v$ of $k$.  If
\begin{equation}\label{com16}
\sum_v t_v \log |\alpha|_v = 0
\end{equation}
for all $\alpha$ in $k^{\times}/\Tor\bigl(k^{\times}\bigr)$, then the function 
$v\rightarrow t_v$ is constant.
\end{lemma}

\begin{proof}  Let $S$ be a finite set of places of $k$ containing all archimedian places, and
assume that the cardinality of $S$ is $s \ge 2$.  We write $\R^s$ for the $s$-dimensional 
real vector space of column vectors $\bx = (x_v)$ having rows indexed by places $v$ in $S$.  
In particular, we write $\bt = (t_v)$ for the column vector in 
$\R^s$ formed from the values of the function $v\rightarrow t_v$ restricted to $S$.  And
we write $\bu = (u_v)$ for the column vector in $\R^s$ such that $u_v = 1$
for each $v$ in $S$.  

Let
\begin{equation*}\label{com17}
U_S(k) = \{\eta \in k: |\eta|_v = 1\ \text{for all}\ v\notin S\}
\end{equation*}
denote the multiplicative group of $S$-units in $k$.  By the $S$-unit theorem (stated 
as \cite[Theorem 3.5]{narkiewicz1974}), there exist multiplicatively independent
elements 
$\xi_1, \xi_2, \dots , \xi_{s-1}$ in $U_S(k)$ which form a fundamental system of
$S$-units.  Write
\begin{equation*}\label{com18}
M = \bigl([k_v:\Q_v] \log \|\xi_r\|_v\bigr)
\end{equation*}
for the associated $(s-1)\times s$ real matrix, where $r = 1, 2, \dots , s-1$ indexes rows and 
$v$ in $S$ indexes columns.  As the $S$-regulator does not vanish, the matrix $M$ has rank 
$(s-1)$.  Hence the null space
\begin{equation*}
\eN = \big\{\bx\in \R^s: M\bx = \bo\big\}
\end{equation*}
has dimension $1$.  From the product formula we have $M \bu = \bo$.  Therefore $\eN$ is
spanned by the vector $\bu$.  By hypothesis we have $M \bt = \bo$, and
it follows that $\bt$ is a scalar multiple of $\bu$.  That is, the function
$v\rightarrow t_v$ is constant on $S$.  As $S$ is arbitrary the lemma is proved.
\end{proof}

We now prove Theorem \ref{thm1}.  Let $\E_1$ denote the closure of $\F$ in $\X$.  As $\F$ is a vector space over
the field $\Q$, it follows that $\E_1$ is a vector space over $\R$, and therefore $\E_1$ is a closed 
linear subspace of $\X$.  If $\E_1$ is a proper subspace then it follows from the 
Hahn-Banach theorem (see \cite[Theorem 3.5]{rudin1991}) that there exists a continuous linear 
functional $\Phi:\X\rightarrow \R$ such that $\Phi$ vanishes on $\E_1$, but $\Phi$ is not the zero linear 
functional on $\X$.  We will show that such a $\Phi$ does not exist, and therefore we must have $\E_1 = \X$.

Let $\Phi:\X\rightarrow \R$ be a continuous linear functional that vanishes on $\E_1$,
but $\Phi$ is not the zero linear functional on $\X$.
It follows from (\ref{space1}) that $\X^{\perp} \subseteq L^{\infty}(Y,\B,\lambda)$ is the 
one dimensional subspace spanned by the constant function $1$.  As the dual space $\X^*$ can
be identified with the quotient space $L^{\infty}(Y,\B,\lambda)/\X^{\perp}$, there exists a
function $\vphi(y)$ in $L^{\infty}(Y,\B,\lambda)$ such that $\vphi(y)$ and the constant function 
$1$ are linearly independent, and
\begin{equation*}
\Phi(F) = \int_Y F(y) \vphi(y)\ \dla(y)
\end{equation*}
for all $F$ in $\X$.  Because $\Phi$ vanishes on $\E_1$ we have
\begin{equation}\label{com30}
\int_Y f_{\alpha}(y) \vphi(y)\ \dla(y) = 0
\end{equation}
for each function $f_{\alpha}$ in $\F$.  

Now let $k$ be a number field in $\eL$ and let $\alpha$ be in 
$k^{\times}/\Tor\bigl(k^{\times}\bigr) \subseteq \G$.  From (\ref{rr4}) and (\ref{com30}) we find that
\begin{align}\label{com31}
\begin{split}
0 &= \sum_v \bigg\{\int_{\psi_k^{-1}(v)} \log \|\alpha\|_y \vphi(y)\ \dla(y)\bigg\} \\
  &= \sum_v \bigg\{\int_{\psi_k^{-1}(v)} \vphi(y)\ \dla(y)\bigg\} \log \|\alpha\|_v \\
  &= \sum_v \bigg\{\lambda\bigl(\psi_k^{-1}(v)\bigr)^{-1} \int_{\psi_k^{-1}(v)} \vphi(y)\ \dla(y)\bigg\} \log |\alpha|_v.
\end{split}
\end{align}
It follows from Lemma \ref{lem6} that the function
\begin{equation*}\label{com32}
v\rightarrow \lambda\bigl(\psi_k^{-1}(v)\bigr)^{-1} \int_{\psi_k^{-1}(v)} \vphi(y)\ \dla(y)
\end{equation*}
is constant on the set of places $v$ of $k$.  We write $c(k)$ for this constant.

Let $k\subseteq l$ be number fields in $\eL$, and let $v$ be a place of $k$.  Using (\ref{rr4}) and (\ref{rr9}) we have
\begin{equation*}\label{com33}
\lambda\bigl(\psi_k^{-1}(v)\bigr) = |W_v(l/k)| \lambda\bigl(\psi_l^{-1}(w)\bigr)
\end{equation*}
for all places $w$ in the set $W_v(l/k)$.  This leads to the identity
\begin{align}\label{com34}
\begin{split}
c(l) &= |W_v(l/k)|^{-1} \sum_{w\in W_v(l/k)} 
		\bigg\{\lambda\bigl(\psi_l^{-1}(w)\bigr)^{-1} \int_{\psi_l^{-1}(w)} \vphi(y)\ \dla(y)\bigg\} \\
     &= \lambda\bigl(\psi_k^{-1}(v)\bigr)^{-1} \sum_{w\in W_v(l/k)} \bigg\{\int_{\psi_l^{-1}(w)} \vphi(y)\ \dla(y)\bigg\} \\
     &= \lambda\bigl(\psi_k^{-1}(v)\bigr)^{-1} \int_{\psi_k^{-1}(v)} \vphi(y)\ \dla(y) \\
     &= c(k).
\end{split}
\end{align}
Thus there exists a real number $C$ such that $C = c(k)$ for all fields $k$ in $\eL$.  

Let $g$ belong to $LC_c(Y)$.  By Lemma \ref{lem4} there exists a number field $l$ in $\eL$ such that $g$ is constant
on $\psi_l^{-1}(w)$ for each place $w$ of $l$.  Therefore we have
\begin{align}\label{com35}
\begin{split}
\int_Y g(y) \vphi(y)\ \dla(y) &= \sum_w \bigg\{\int_{\psi_l^{-1}(w)} g(y) \vphi(y)\ \dla(y)\bigg\} \\
     &= C \sum_w \Big\{\lambda\bigl(\psi_l^{-1}(w)\bigr) g\bigl(\psi_l^{-1}(w)\bigr)\Big\} \\
     &= C \sum_w \bigg\{\int_{\psi_l^{-1}(w)} g(y)\ \dla(y)\bigg\} \\
     &= C \int_Y g(y)\ \dla(y).
\end{split}
\end{align}
By Lemma \ref{lem5} the set $LC_c(Y)$ is dense in $L^1(Y,\B,\lambda)$, and we conclude from (\ref{com35}) that
\begin{equation*}\label{com36}
\int_Y F(y) \vphi(y)\ \dla(y) = C \int_Y F(y)\ \dla(y)
\end{equation*}
for all $F$ in $L^1(Y,\B,\lambda)$.  This shows that $\vphi(y) = C$ in $L^{\infty}(Y,\B,\lambda)$, and
so contradicts our assumption that $\vphi(y)$ and the constant function $1$ are linearly independent.  Hence the
continuous linear functional $\Phi$ does not exist, and therefore $\E_1 = \X$.

%%%%%%%%%%%%%%%%%%%%%%%%%%%%%%%%%%%%%%%%%%%%%%%%%%%%%%%%%%%%%%%%%%%%%%%%%%%%%%%%%%
%%%%%%%%%%%%%%%%%%%%%%%%%%%%%%%%%%%%%%%%%%%%%%%%%%%%%%%%%%%%%%%%%%%%%%%%%%%%%%%%%%

\section{Proof of Theorem \ref{thm2} and Theorem \ref{thm3}}

We suppose that $1 < p < \infty$ and write $\E_p$ for the closure of $\F$ in $L^p(Y,\B,\lambda)$.
As before, $\E_p$ is a closed linear subspace.  By the Hahn-Banach theorem it suffices to show that if
$\Phi:L^p(Y,\B,\lambda)\rightarrow \R$ is a continuous linear functional that vanishes on $\E_p$, then
in fact $\Phi$ is identically zero on $L^p(Y,\B,\lambda)$.

Let $p^{-1} + q^{-1} = 1$, and let $\vphi(y)$ be an element of $L^q(Y,\B,\lambda)$ such that
\begin {equation*}\label{com40}
\Phi(F) = \int_Y F(y) \vphi(y)\ \dla(y)
\end{equation*}
for all $F$ in $L^p(Y,\B,\lambda)$.  We assume that $\Phi$ vanishes on $\E_p$, and then we have
\begin{equation}\label{com41}
\int_Y f_{\alpha}(y) \vphi(y)\ \dla(y) = 0
\end{equation}
for each function $f_{\alpha}$ in $\F$.  

Let $k$ be a number field in $\eL$ and let $\alpha$ be in $k^{\times}/\Tor\bigl(k^{\times}\bigr) \subseteq \G$. 
As before we apply (\ref{rr4}) and (\ref{com30}) to obtain the identity (\ref{com31}). Then Lemma \ref{lem6} implies
that the function
\begin{equation}\label{com42}
v\rightarrow \lambda\bigl(\psi_k^{-1}(v)\bigr)^{-1} \int_{\psi_k^{-1}(v)} \vphi(y)\ \dla(y)
\end{equation}
is constant on the set of places $v$ of $k$.  Now, however, we apply H\"older's inequality and find that
\begin{align*}\label{com43}
\begin{split}
\sum_v~\biggl|\lambda\bigl(\psi_k^{-1}(v)\bigr)^{-1}&\int_{\psi_k^{-1}(v)} \vphi(y)\ \dla(y)\biggr|^q \\
        &\le \sum_v~\bigg\{\lambda\bigl(\psi_k^{-1}(v)\bigr)^{-1} \int_{\psi_k^{-1}(v)} |\vphi(y)|^q\ \dla(y)\bigg\} \\
        &\le [k:\Q] \int_Y |\vphi(x)|^q\ \dla(y) < \infty.
\end{split}
\end{align*}
This shows that the constant value of the function (\ref{com42}) is zero.  Thus we have
\begin{equation*}\label{com44}
\int_{\psi_k^{-1}(v)} \vphi(y)\ \dla(y) = 0
\end{equation*}
for all $k$ in $\eL$ and for all places $v$ of $k$.  It follows using Lemma \ref{lem4} that
\begin{equation*}\label{com45}
\int_Y g(y) \vphi(y)\ \dla(y) = 0
\end{equation*}
for all $g$ in $LC_c(Y)$.  By Lemma \ref{lem5} the set $LC_c(Y)$ is dense in $L^p(Y,\B,\lambda)$, and we conclude that
the continuous linear functional $\Phi$ is identically zero.  This completes the proof of Theorem \ref{thm2}.

Next we suppose that $\E_{\infty}$ is the closure of $\F$ in $C_0(Y)$.  Again it suffices to show that if 
$\Phi:C_0(Y)\rightarrow \R$ is a continuous linear functional that vanishes on $\E_{\infty}$, then
$\Phi$ is identically zero on $C_0(Y)$.  If $\Phi$ is such a linear functional, then by the Riesz representation
theorem (see \cite[Theorem 6.19]{rudin1987}) there exists a regular signed measure $\nu$, defined on the $\sigma$-algebra
$\B$ of Borel sets in $Y$, such that
\begin{equation*}\label{com50}
\Phi(F) = \int_Y F(y)\ \dnu(y)
\end{equation*}
for all $F$ in $C_0(Y)$.  Moreover, we have $\|\Phi\| = \|\nu\|$, where $\|\Phi\|$ is the norm of the linear functional
$\Phi$ and $\|\nu\|$ is the total variation of the signed measure $\nu$. 
We assume that $\Phi$ vanishes on $\E_{\infty}$, and therefore
\begin{equation*}\label{com51}
\int_Y f_{\alpha}(y)\ \dnu(y) = 0
\end{equation*}
for each function $f_{\alpha}$ in $\F$.  By arguing as in the proof of Theorem 2, we conclude that for each number field
$k$ in $\eL$ the function
\begin{equation}\label{com52}
v\rightarrow \lambda\bigl(\psi_k^{-1}(v)\bigr)^{-1} \nu\bigl(\psi_k^{-1}(v)\bigr), 
\end{equation}
defined on the set of all places $v$ of $k$, is constant.  As
\begin{align*}\label{com53}
\begin{split}
\sum_v~\Bigl|\lambda\bigl(\psi_k^{-1}(v)\bigr)^{-1}&\nu\bigl(\psi_k^{-1}(v)\bigr)\Bigr| \\
    &\le [k:\Q] \sum_v~\bigl|\nu\bigl(\psi_k^{-1}(v)\bigr)\bigr| \\
    &\le [k:\Q] \|\nu\| < \infty,
\end{split}
\end{align*}
we conclude that the value of the constant function (\ref{com52}) is zero.  This shows that
\begin{equation*}\label{com54}
\nu\bigl(\psi_k^{-1}(v)\bigr) = 0
\end{equation*}
for all $k$ in $\eL$ and for all places $v$ of $k$.  It follows as before that
\begin{equation*}\label{com55}
\Phi(g) = \int_Y g(y)\ \dnu(y) = 0
\end{equation*}
for all $g$ in $LC_c(Y)$.  As $LC_c(Y)$ is dense in $C_0(Y)$ by Lemma \ref{lem5}, we find that $\Phi$ is identically
zero on $C_0(Y)$.  This proves Theorem \ref{thm3}.
%%%%%%%%%%%%%%%%%%%%%%%%%%%%%%%%%%%%%%%%%%%%%%%%%%%%%%%%%%%%%%%%%%%%%%%%%%%%%%%%
%%%%%%%%%%%%%%%%%%%%%%%%%%%%%%%%%%%%%%%%%%%%%%%%%%%%%%%%%%%%%%%%%%%%%%%%%%%%%%%%

%\bibliography{disrefs}

\end{document}